\documentclass{article}

\title{The Cantor-Schr\"oder-Bernstein Theorem for $\infty$-groupoids}
\author{Mart\'{\i}n H\"otzel Escard\'o
  \\ {\small School of Computer Science}
  \\ {\small University of Birmingham, UK}}
\usepackage{amsmath,stmaryrd,amssymb,amsthm,url,hyperref,color}

\definecolor{links}{rgb}{0.47, 0.27, 0.23}

\hypersetup{
    colorlinks, linkcolor={links},
    citecolor={links}, urlcolor={links}
}

\swapnumbers
\newtheorem{theorem}{Theorem}[section]
\newtheorem{claim}[theorem]{Claim}

\theoremstyle{definition}

\newtheorem{example}[theorem]{Example}
\newtheorem{remark}[theorem]{Remark}

\begin{document}

\maketitle

\begin{abstract}
  We show that the Cantor-Schr\"oder-Bernstein Theorem for homotopy
  types, or $\infty$-groupoids holds in the following form: For any
  two types, if each one is embedded into the other, then they are
  equivalent. The argument is developed in the language of homotopy
  type theory, or Voevodsky's univalent foundations (HoTT/UF), and
  requires classical logic. It follows that the theorem holds in any
  boolean $\infty$-topos.
\end{abstract}

\section{Introduction}

The classical Cantor-Schr\"oder-Bernstein Theorem of set theory, formulated by Cantor and first proved by Bernstein, states that for any pair of sets, if there is an injection of each one into the other, then the two sets are in bijection.
There are proofs that use excluded middle but not choice. That excluded middle is absolutely necessary was recently established Pierre Pradic and Chad E. Brown~\cite{pradic:brown}.

The appropriate \emph{principle of excluded middle} for HoTT/UF~\cite{hottbook} says that every subsingleton (or proposition, or truth value) is either empty or pointed. The statement that \emph{every type} is either empty or pointed is much stronger, and amounts to \emph{global choice}, which is incompatible with univalence~\cite[Theorem 3.2.2]{hottbook}. In fact, in the presence of global choice, every type is a set by Hedberg's Theorem, but univalence gives types that are not sets. Excluded middle middle, however, is known to be compatible with univalence, and is validated in Voevodsky's model of simplicial sets. And so is (non-global) choice, but it is not needed for our purposes.

Even assuming excluded middle, it may seem unlikely at first sight
that the Cantor-Schr\"oder-Bernstein Theorem (CSB) can be generalized from sets
to arbitrary homotopy types, or $\infty$-groupoids:
\begin{enumerate}
 	\item CSB fails for 1-categories.
In fact, it already fails for posets. For example, the intervals $(0,1)$ and $[0,1]$ are order-embedded into each other, but they are not order isomorphic, or equivalent as categories.
 	\item The known proofs of CSB for sets rely on deciding equality of elements of sets, but, in the presence of excluded middle, the types that have decidable equality are precisely the sets, by Hedberg's Theorem.
\end{enumerate}
In set theory, a map $f : X \to Y$ is an injection if and only if it is left-cancellable, in the sense that $f(x)=f(x')$ implies $x=x'$. But, for types $X$ and $Y$ that are not sets, this notion is too weak, and, moreover, is not a proposition as the identity type $x = x'$ has multiple elements in general. The appropriate notion of \emph{embedding} for a function $f$ of \emph{arbitrary} types $X$ and $Y$ is given by any of the following two equivalent conditions:
\begin{enumerate}
 	\item The map $\operatorname{ap}(f,x,x') : x = x' \to f(x) = f(x')$ is an equivalence for any $x,x':X$.
 	\item The fibers of $f$ are all subsingletons.
\end{enumerate}

A map of sets is an embedding if and only if it is left-cancellable. However, for example, any map $1 \to Y$ that picks a point $y:Y$ is left-cancellable, but it is an embedding if and only if the point $y$ is homotopy isolated, which amounts to saying that the identity type $y = y$ is contractible. This fails, for instance, when the type $Y$ is the homotopical circle $S^1$, for any point $y$, or when $Y$ is a univalent universe and $y:Y$ is the two-point type, or any type with more than one automorphism.

\begin{example}[Pradic~\cite{pradic:example}]
  There is a pair of left-cancellable maps between the types
  $\mathbb{N} \times S^1$ and $1 + \mathbb{N} \times S^1$ (taking
  $\operatorname{inl}$ going forward and, going backward, mapping
  $\operatorname{inl}(*)$ to $(0, \operatorname{base})$ and shifting
  the indices of the circles by one), but no equivalence between these
  two types.
\end{example}

\section{Cantor-Schr\"oder-Bernstein for $\infty$-groupoids}

As explained in the introduction, our argument is in the language of
HoTT/UF and requires classical logic. Because HoTT/UF can be
interpreted in any $\infty$-topos~\cite{shulman:infty}, it follows
that the following theorem holds in any \emph{boolean} $\infty$-topos.
We assume the terminology and notation of the HoTT
book~\cite{hottbook}.

\begin{theorem}
For any two types, if each one is embedded into the
other, then they are equivalent, in the presence of excluded middle.
\end{theorem}

We adapt Halmos' proof~\cite{MR0453532} for sets. We need to
refomulate the argument so that excluded middle is applied to
truth-valued, rather than type-valued, mathematical statements, and
this is the contribution in this note (see Remark~\ref{main:remark}
below).  We don't need to invoke univalence, the existence of
propositional truncations or any other higher inductive type for our
construction. But we do rely on function extensionality. An
Agda~\cite{agda} version of the following argument is
available~\cite{2019arXiv191100580H,escardo:csb}.

\begin{proof}
Let $f : X \to Y$ and $g : Y \to X$ be embeddings of arbitrary types $X$ and~$Y$.
We say that $x:X$ is a $g$-point if for any $x_0 : X$ and
$n : \mathbb{N}$ with $(g \circ f)^n (x_0)=x$, the $g$-fiber of $x_0$
is inhabited. Using the assumption that $g$ is an embedding, we see
that being a $g$-point is property rather than data, because
subsingletons are closed under products by function extensionality.

Considering $x_0=x$ and $n=0$, we see that if $x$ is a $g$-point then
the $g$-fiber of $x$ is inhabited, and hence we get a function
$g^{-1}$ of $g$-points of $X$ into $Y$. By construction, we have that
$g(g^{-1}(x))=x$. In particular, if $g(y)$ is a $g$-point for a given
$y:Y$, we conclude that $g(g^{-1}(g(y)))=g(y)$, and because $g$, being
an embedding, is left-cancellable, we get $g^{-1}(g(y))=y$.

Now define $h:X \to Y$ by
\[
  h(x) = \begin{cases}
           g^{-1}(x) & \text{if $x$ is a $g$-point}, \\
           f(x) & \text{otherwise.}
         \end{cases}
\]
To conclude the proof, it is enough to show that $h$ is left-cancellable and split-surjective, as any such map is an equivalence.

To see that $h$ is left-cancellable, it is enough to show that the
images of $f$ and $g^{-1}$ in the definition of $h$ are disjoint,
because $f$ and $g^{-1}$ are left-cancellable. For that purpose, let
$x$ be a non-$g$-point and $x'$ be a $g$-point, and, for the sake of
contradiction, assume $f(x) = g^{-1}(x')$. Then
$g(f(x))=g(g^{-1}(x'))=x'$. Now, because if $g(f(x))$ were a $g$-point
then so would be $x$, we conclude that it isn't, and hence neither is
$x'$, which contradicts the assumption.

To see that $h$ is a split surjection, say that $x : X$ is an $f$-point if there are designated $x_0 : X$ and $n : \mathbb{N}$ with $(g \circ f)^n (x_0)=x$ and the $g$-fiber of $x_0$ empty. This is data rather than property, and so this notion could not have been used for the construction of $h$. But every non-$f$-point is a $g$-point, applying excluded middle to the $g$-fiber of $x_0$ in the definition of $g$-point.

\begin{claim}
  If $g(y)$ is not a $g$-point, then there is a designated point
  $(x,p)$ of the $f$-fiber of $y$, with $x : X$ and $p : f(x)=y$,
  such that $x$ is not a $g$-point either.
\end{claim}

To prove the claim, first notice that it is impossible that $g(y)$ is not an $f$-point, by the above observation. But this is not enough to conclude that it is an $f$-point, because excluded middle applies to subsingletons only, which the notion of $f$-point isn't. However, it is readily seen that if $g(y)$ is an $f$-point, then there is a designated point $(x,p)$ in the $f$-fiber of $y$. From this it follows that it impossible that the subtype of the fiber consisting of the elements $(x,p)$ with $x$ not a $g$-point is empty. But the $f$-fiber of $y$ is a proposition because $f$ is an embedding, and hence so is the subtype, and therefore the claim follows by double-negation elimination.

We can now resume the proof that $h$ is a split surjection. For any $y:Y$, we check whether $g(y)$ is a $g$-point. If it is, we map $y$ to $g(y)$, and if it isn't we map $y$ to the point $x : X$ given by the claim, which concludes the proof of the theorem.
\end{proof}

\begin{remark} \label{main:remark} So, in this argument we don't apply
  excluded middle to equality directly, which we wouldn't be able to
  as the types $X$ and $Y$ are not necessarily sets. We instead apply
  it to (1) the property of being a $g$-point, defined in terms of the
  fibers of $g$, to define $h$, (2) a fiber of $g$, and (3) a subtype
  of a fiber of $f$. These three types are propositions because the
  functions $f$ and $g$ are embeddings rather than merely
  left-cancellable maps.
\end{remark}

\begin{remark}
  If the type $X$ in the proof is connected, then every map of $X$
  into a set is constant. In particular, the property of being a
  $g$-point is constant, because the type of truth values is a set
  (assuming univalence for subsingletons). Hence, by excluded middle,
  it is constantly true or constantly false, and so $h=g^{-1}$ or
  $h = f$, which means that one of the embeddings $f$ and $g$ is
  already an equivalence.  Mike Shulman (personal communication)
  observed that this is true even without excluded middle: If $X$ is
  connected and we have an embedding $g : Y \to X$ and any function at
  all $f : X \to Y$, then $g$ is an equivalence. For any $x:X$, we
  have $\left\lVert g(f(x)) = x \right\rVert$ since $X$ is connected;
  thus $g$ is (non-split) surjective. But a surjective embedding is an
  equivalence.
\end{remark}

\bibliographystyle{plain}
\bibliography{csb}

\vfill

% \noindent Mart\'{\i}n H\"otzel Escard\'o \\
% \noindent School of Computer Science \\
% \noindent University of Birmingham, UK

\end{document}